\documentclass[12pt,pdflatex,sn-mathphys-num]{sn-jnl}
\usepackage{graphicx}%
\usepackage{multirow}%
\usepackage{amsmath,amssymb,amsfonts}%
\usepackage{enumitem}%
\setlist[itemize]{leftmargin=2.5em, topsep=0pt, partopsep=0pt, parsep=0pt}
\usepackage{amsthm}%
\usepackage{mathrsfs}%
\usepackage[title]{appendix}%
\usepackage{textcomp}%
\usepackage{manyfoot}%
\usepackage{booktabs}%
\usepackage[table]{xcolor}
\usepackage{algorithm}%
\usepackage{algorithmicx}%
\usepackage{algpseudocode}%
\usepackage{listings}%
\usepackage{caption} 
\usepackage[labelformat=simple]{subcaption}
\usepackage{colortbl}
\usepackage{float} 
\usepackage[section]{placeins} 
\usepackage{setspace}
\usepackage{lineno}

\theoremstyle{thmstyleone}%
\theoremstyle{thmstyletwo}%

\theoremstyle{thmstylethree}%

\graphicspath{{./FIGURES/}} 

\begin{document}

\title[Nonlinear Spectral Modeling and Control]{%
  Nonlinear Spectral Modeling and Control of
  Soft-Robotic Muscles from Data%
}


\author[1]{\fnm{Leonardo} \sur{Bettini}}\equalcont{These authors contributed equally to this work.}

\author[2]{\fnm{Amirhossein} \sur{Kazemipour}}\equalcont{These authors contributed equally to this work.}

\author*[2]{\fnm{Robert K.} \sur{Katzschmann}}

\author*[1]{\fnm{George} \sur{Haller}}

\email{rkk@ethz.ch; georgehaller@ethz.ch}

\affil[1]{Institute for Mechanical Systems, ETH Zurich, Leonhardstrasse 21, 8092 Zurich, Switzerland}

\affil[2]{Soft Robotics Lab, ETH Zurich, Tannenstrasse 3, 8092 Zurich, Switzerland}

\abstract{Artificial muscles are essential for compliant musculoskeletal robotics but complicate control due to nonlinear multiphysics dynamics. Hydraulically amplified electrostatic (HASEL) actuators, a class of soft artificial muscles, offer high performance but exhibit memory effects and hysteresis. Here we present a data-driven reduction and control strategy grounded in spectral submanifold (SSM) theory. In the adiabatic regime, where inputs vary slowly relative to intrinsic transients, trajectories rapidly converge to a low-dimensional slow manifold. We learn an explicit input-to-output map on this manifold from forced-response trajectories alone, avoiding decay experiments that can trigger hysteresis. We deploy the SSM-based model for real-time control of an antagonistic HASEL-clutch joint. This approach yields a substantial reduction in tracking error compared to feedback-only and feedforward-only baselines under identical settings. This record-and-control workflow enables rapid characterization and high-performance control of soft muscles and muscle-driven joints without detailed physics-based modeling.}




\maketitle

\section{Introduction}

The adaptability, compliance, and efficiency of biological systems increasingly inspire advances in robotics. In musculoskeletal systems like the human arm, rigid elements (bones) provide structure and transmit forces, while soft components (muscles, ligaments, tendons) connect and actuate them, enabling fine control and safe interaction with delicate environments. Soft elements prevent damage from rigid contact, ensuring compliance. Soft-material-based robots extend these advantages, offering flexibility, adaptability, and compatibility with complex geometries, particularly valuable in human-interactive applications such as surgical tools and prosthetics. Artificial muscles, emulating biological force generation, are key to innovations in musculoskeletal robotics and assistive technologies.

Despite these benefits, modeling and control of soft artificial muscles remain challenging, primarily due to the complex interplay among soft materials, electrostatic forces, fluidic dynamics, and frictional contacts. These effects lead to memory, hysteresis, and other nonlinearities that complicate analysis and real-time control. Reduced-order models (ROMs) are therefore attractive, as they capture essential dynamics with minimal equations while balancing interpretability and efficiency (see \citep{benner2015,ghadami2022} for recent reviews).
In practice, artificial muscles operate within larger assemblies and are affected by fabrication imperfections such as partial adhesion or misalignment, which add uncertainty. Under these conditions, analytic ROMs become limited and hard to scale, motivating a shift toward data-driven modeling that captures system dynamics directly from data without precise parameter identification.

\begin{figure*}[!htbp]
    \centering
    \includegraphics[scale=1]{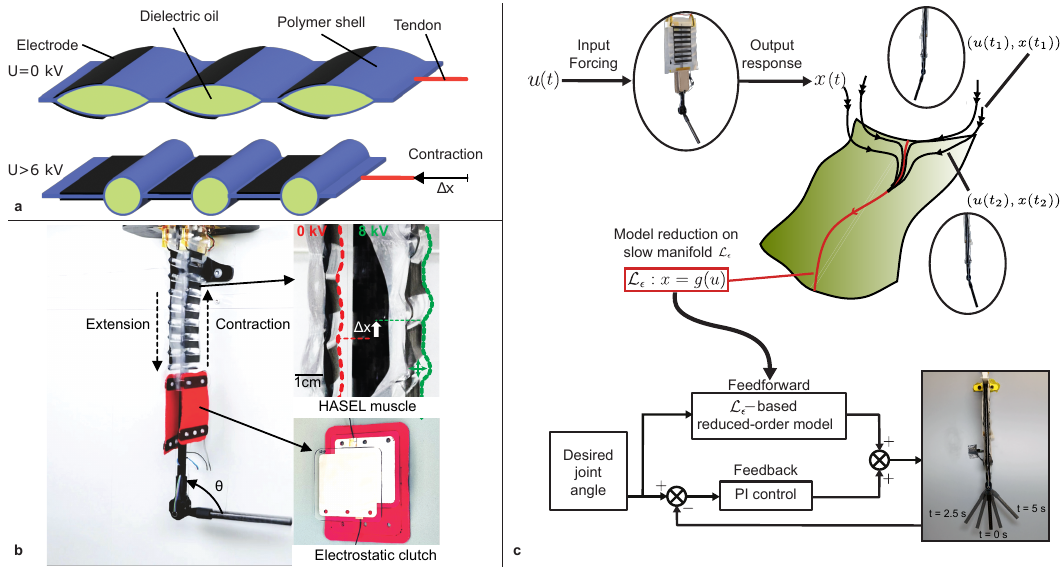}
     \caption{\textbf{Overview of data-driven slow manifold modeling and control for antagonistic artificial muscles.} 
(\textbf{a}) HASEL actuator working principle: a flexible polymer shell filled with dielectric oil and covered by electrodes; applying voltage (up to $U=8$~kV) generates electrostatic pressure that redistributes the fluid, producing axial contraction $\Delta x$ against a tendon load. 
(\textbf{b}) Experimental platform: antagonistic musculoskeletal joint actuated by HASEL artificial muscles paired with electrostatic clutches, enabling bidirectional motion and variable stiffness; $\theta$ denotes the joint angle. The inset shows the contraction $\Delta x$ at 0~kV and 8~kV. 
(\textbf{c}) Data-driven modeling and control: pronounced time-scale separation allows trajectories to rapidly converge onto a low-dimensional slow manifold $\mathcal{L}_\epsilon$, learned directly from forced-response trajectories. The learned polynomial map $x = g(u)$ relates control input $u(t)$ to system state $x(t)$. The inverse slow-manifold model provides feedforward control, combined with proportional-integral (PI) feedback for disturbance rejection.} 
    \label{method_preview}
\end{figure*}

Among data-driven modeling techniques, the recently developed spectral submanifold (SSM) reduction has addressed similar challenges across various physical systems \citep{haller_2025}. SSM reduction preserves the essential nonlinear dynamics of high-dimensional systems in a compact and interpretable form. Specifically, slow SSMs are low-dimensional, smooth, and attracting invariant manifolds in the phase space, emerging from the system’s slowest linear modes \citep{haller2016}. Trajectories near the SSM quickly converge onto it and evolve according to its internal dynamics, yielding accurate and robust nonlinear reduced-order models. SSM reduction has proven effective in both equation-based and purely data-driven settings \citep{haller_2025}, accurately capturing nonlinear decay in free responses and predicting forced responses under periodic and quasiperiodic excitations \citep{cenedese2022,cenedese2022_mechanical,joar_2023,bettini_2024,joar_2024,bettini_2025}.

Recent theoretical advances have extended SSM theory to nonlinear systems under temporally aperiodic excitations \citep{haller_2024,haller_2025}, broadening its applicability to practical engineering fields such as robotics. The extended theory accommodates time-dependent inputs that are either moderate in amplitude or vary slowly relative to the system’s internal time scales of the uncontrolled system. The latter setting case is known as the adiabatic setting. Building on promising simulation results in soft-robot control via SSM theory \citep{Alora_2023}, \citet{Alora_2025} demonstrated an experimental implementation using a leading-order approximation of an adiabatic SSM (aSSM) anchored at the origin, achieving both computational and data efficiency. However, this approach assumes trajectories remain near the origin, limiting its scope. This restriction was recently lifted by \citet{kaundinya_2025}, who developed a data-driven aSSM reduction and control strategy for a finite-element model of a soft trunk robot.

In this work, we also focus on the adiabatic setting, where intrinsic dynamics evolve much faster than external forcing. We further consider the regime in which this time-scale separation is strong enough to justify neglecting the internal dynamics on the adiabatic spectral submanifold (aSSM), focusing instead on the slow manifold (SM) to which the aSSM is anchored. This is a fair approximation for soft artificial muscles, whose internal actuation and material transients are typically strongly damped and decay rapidly relative to the motion along a desired trajectory. Crucially, the learned SM admits a fast inverse map, enabling real-time feedforward control that can be combined with feedback for accurate trajectory tracking in antagonistic muscle configurations.

Experimental evidence of pronounced time-scale separation in soft-robotic systems comes, for example, from Hydraulically Amplified Self-Healing Electrostatic (HASEL) actuators, which are also the focus of this work. HASELs are part of the broader class of electrohydraulic artificial muscles, offering a compelling alternative to motor-driven systems by overcoming the bulkiness, rigidity, and limited compliance of electromagnetic motors \citep{Asano_2016,Jantsch_2013,Richter_2016,acome_2018,Yang_2019,Buchner_2024,Gravert_2024}.
These designs are inspired by animal musculoskeletal architecture, where bones, tendons, and contractile actuators work together for efficient, compliant motion \citep{Kurumaya_2016,Niiyama_2007,Niiyama_2010,Ogawa_2011}. HASEL actuators combine soft fluidic and electrostatic actuation to achieve compliant yet powerful contraction in lightweight, scalable structures. They consist of oil-filled polymer pouches partially coated with electrodes; applying voltage deforms the pouch, redistributing the dielectric liquid and generating contraction. Extension occurs passively under external loads or via antagonistic actuation in multi-muscle assemblies \citep{Kazemipour_2024}, enabling bidirectional motion, variable stiffness, and rapid disturbance rejection. The fast electrostatic and fluidic transients relative to the slower mechanical response produce the slow–fast dynamics of the adiabatic regime, making HASELs an ideal platform for the SM-based modeling introduced here.

\citet{Kellaris_2019} proposed a quasi-static model for HASEL actuators, wherein the actuator stroke, parametrized by a single degree of freedom, is obtained by minimizing the system’s free energy to find equilibrium. Originally developed for electrostatic zipping actuators \citep{Saif_1999,Maffli_2013,Righi_2018}, this approach links static deformation to applied voltage through an analytical relationship depending on physical parameters like initial pouch angle, film thickness, actuator width, undeformed pouch length, and dielectric permittivity, requiring precise measurements for accuracy. Later studies highlighted modeling’s role in closed-loop control. \citet{Johnson_2020} identified the input-output relationship via frequency response testing, enabling frequency-domain controller design. \citet{Yeh_2022} developed a three-dimensional dynamical model using the Port-Hamiltonian formalism. Like the quasi-static model, these approaches depend on detailed physical parameter identification, illustrating the trade-off between interpretability and the practical challenges of measuring precise system properties.

These studies, however, focus on isolated HASEL actuators, which are rarely used alone. In practice, HASELs are embedded in complex assemblies. For example, \citet{Kazemipour_2024} enhanced a robotic arm’s range of motion by integrating HASELs with electrostatic clutches. These thin, flexible electroadhesive elements modulate shear adhesion via a low-current voltage signal, rapidly switching between a high-holding-force state (to transmit muscle force) and a low-friction state (to release the opposing tendon), thereby enabling bidirectional joint motion and variable stiffness without tendon slack. Precise control of such antagonistic systems requires accurate, computationally efficient models. Continuum models are too slow for real-time use, while black-box approaches lack reliability. Multi-physics coupling and switching dynamics introduce nonlinearities and hysteresis that resist analytical treatment, motivating data-driven reduced-order modeling: learning system dynamics directly from experiments while remaining computationally efficient. In this context, our SM-based, data-driven method provides a scalable, interpretable approach for predicting and controlling multi-actuator HASEL systems, inherently accounting for component interactions and fabrication imperfections.

In this work, we introduce a methodology for performing SM reduction directly from experimental data and demonstrate its use for real-time control of antagonistic artificial muscles. As outlined in Fig. \ref{method_preview}c, the procedure is detailed in Sections \ref{SSM_section} and \ref{aSSM_section}, applied to numerical examples in Section \ref{analytical_models}, and validated on experimental data in Section \ref{data_driven_modeling}. Section \ref{control_section} presents closed-loop control experiments on an antagonistic musculoskeletal joint with HASEL actuators and electrostatic clutches. We compare three strategies, PI feedback only, inverse-SSM feedforward only, and their combination (Methods \ref{methods_control}), showing that the combined controller achieves substantially lower tracking errors. On randomized test trajectories under identical saturation and slew-rate limits, the inverse-SSM+PI controller (RMS error of 2.38$^\circ$, max 10.41$^\circ$)  outperforms the baselines, i.e., feedforward-only (RMS error of 3.63$^\circ$, max 16.67$^\circ$) and feedback-only (RMS error of 7.63$^\circ$, max 26.00$^\circ$). These results confirm that data-driven slow-manifold models not only predict system behavior accurately but also enable high-performance, real-time control of complex antagonistic muscle systems, bridging the gap between theoretical model reduction and practical robotic applications. To the best of our knowledge, this represents the first experimental demonstration of spectral submanifold theory applied to closed-loop control of antagonistic soft robotic systems with quantified performance improvements over standard control baselines.

\section{Results}\label{results}

We first present the theoretical and data-driven reduction method (SSM, aSSM, and SM), validate it on analytic models and on experimental HASEL actuators, and then demonstrate closed-loop control on an antagonistic joint. Together, these results show that a simple inverse map learned on the slow manifold yields accurate predictions and real-time control on hardware.


\subsection{Model reduction to spectral submanifolds}\label{SSM_section}
Consider an autonomous dynamical system of the form 
\begin{equation}\label{autonomous_DS}
\begin{aligned}
&
\boldsymbol{\dot{x}} =\boldsymbol{f}\left(\boldsymbol{x}\right) = \boldsymbol{A}\boldsymbol{x}  + \boldsymbol{f}_0\left(\boldsymbol{x} \right),& \boldsymbol{x} \in \mathbb{R}^n, &\quad \boldsymbol{A} \in \mathbb{R}^{n \times n}, \\&& \boldsymbol{f}_0 \in \mathcal{C}^r, &\quad \boldsymbol{f}_0 = \mathcal{O}\left(\left|\boldsymbol{x} \right|^2 \right),
\end{aligned}
\end{equation}
for some integer $r \ge 1$, where $\boldsymbol{f}_0\left( \boldsymbol{x}\right)$ represents the nonlinear part and $\boldsymbol{f}_0\left(0\right) = \boldsymbol{0}$, so that the system has a fixed point at the origin. We assume the fixed point to be linearly asymptotically stable, i.e.,
\begin{equation}\label{eigenvalues_order}
Re \lambda_n \leq Re \lambda_{n-1} \leq \dots \leq Re \lambda_2 \leq \lambda_1 < 0.
\end{equation}
Each distinct eigenvalue $\lambda_j$ of the linear part $\boldsymbol{A}$ is associated with a real \textit{eigenspace} $E_j$, which is spanned by the real and imaginary parts of the eigenvector or generalized eigenvectors corresponding to $\lambda_j$. A slow \textit{spectral subspace} $E$ is defined as the direct sum of a group of $\ell$ eigenspaces corresponding to the $\ell$ slowest eigenvalues in eq. \eqref{eigenvalues_order}
\begin{equation}
E = E_1 \oplus E_2 \oplus \dots \oplus E_\ell.
\end{equation} 

As discussed in \citet{haller2016}, eigenspaces and spectral subspaces admit unique smoothest nonlinear continuations under generically satisfied nonresonance conditions of the spectrum of $\boldsymbol{A}$. These continuations are referred to as \textit{primary} SSMs and denoted as $\mathcal{W}(E)$. The primary SSM attracts all nearby trajectories in the phase space exponentially and hence its internal dynamics serve as an ideal nonlinear reduced-order model of the full system.

\subsection{Adiabatic spectral submanifold and slow manifold}\label{aSSM_section}

\begin{figure*}[!htbp]
    \centering
    \includegraphics[scale=1]{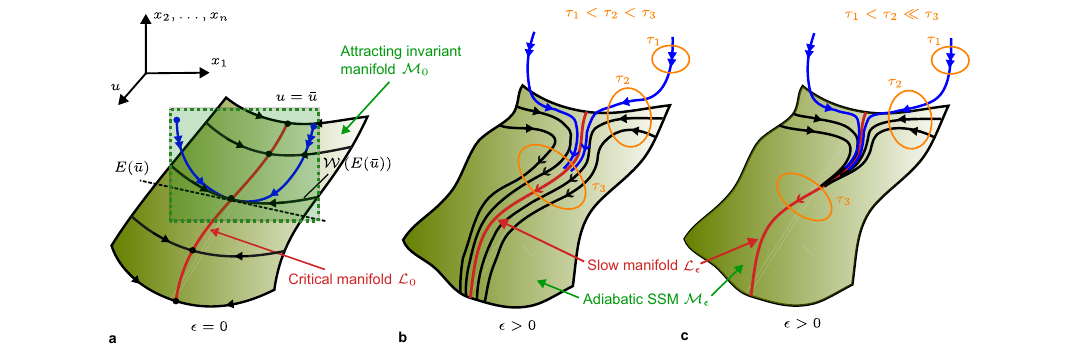}
    \caption{
    (\textbf{a}) Geometry of the attracting invariant manifold $\mathcal{M}_0$ formed by a family of SSMs, $\mathcal{W}(E(u))$, parametrized by the external forcing parameter $u$, for $\epsilon = 0$. For simplicity, $u$ is assumed to be a scalar in this plot. The critical manifold $\mathcal{L}_0$ consists of the collection of the fixed points corresponding to each value of $u$.
    (\textbf{b}) For $\epsilon \ne 0$, the manifold $\mathcal{M}_0$ perturbs into an aSSM, $\mathcal{M}_\epsilon$, and $\mathcal{L}_0$ perturbs into a slow manifold, $\mathcal{L}_\epsilon$. The dynamics of the external forcing is associated with the hollow arrowhead.
    (\textbf{c}) Trajectories rapidly converge to $\mathcal{L}_\epsilon$ when the timescales $\tau_1$ and $\tau_2$ of the autonomous dynamics are much faster than the timescale $\tau_3$ of the external excitation.}
    \label{adiabatic_cases}
\end{figure*}

Let us now reconsider system \eqref{autonomous_DS} under time-dependent forcing $\boldsymbol{u}(t) \in \mathbb{R}^f$ in the form
\begin{equation}\label{nonautonomous_DS}
\begin{array}{ll}
\dot{\boldsymbol{x}} = \boldsymbol{A}\boldsymbol{x} + \boldsymbol{f}_0\left(\boldsymbol{x} \right) + \boldsymbol{f}_1\left(\boldsymbol{x}, \boldsymbol{u} \right), & \quad \boldsymbol{f}_1 \in \mathcal{C}^r, \\[1.2ex]
\dot{\boldsymbol{u}} = \epsilon \boldsymbol{w}\left(\boldsymbol{x}, \epsilon t \right), & \quad 0 \leq \epsilon \ll 1,
\end{array}
\end{equation}
where $\epsilon$ characterizes the order of the rate of change of the external forcing. For $0 \leq \epsilon \ll 1$, $\boldsymbol{u}(t)$ evolves on a slower timescale than the system’s intrinsic dynamics, placing the system in the \textit{adiabatic} regime, as defined by \citet{haller_2024}. In control applications, the term $\boldsymbol{f}_1$ represents the external input applied to the system, with its additive contribution to the autonomous dynamics is justified in the Supplementary Material. This input depends both on the system state, through a feedback law, and on time, to follow a desired trajectory. The adiabatic assumption is particularly well suited to robotics, where the internal dynamics of the system (e.g., actuator response) are significantly faster than the evolution of the controlled motion. This separation of time scales enables the system to rapidly synchronize with the control input and track complex trajectories with high accuracy and responsiveness.

For $\epsilon = 0$, the external forcing becomes constant and acts as a fixed parameter of the system. As seen in Fig. \ref{adiabatic_cases}a, for each fixed value of $\boldsymbol{u}$, the system dynamics evolve toward a unique attracting fixed point along a corresponding slow SSM, $\mathcal{W}\left(E\left(\bar{u} \right)\right)$. The set of these fixed points forms an invariant manifold denoted by $\mathcal{L}_0$, which is a critical manifold in the terminology of geometric singular perturbation theory (see \citet{fenichel_1979}). The union of SSMs taken over all values of $\boldsymbol{u}$ defines a higher-dimensional invariant manifold $\mathcal{M}_0$, which globally attracts nearby trajectories. 

For $\epsilon > 0$, the external forcing varies slowly in time, as shown in Fig. \ref{adiabatic_cases}b,c. In this regime, the critical manifold $\mathcal{L}_0$ perturbs into a slow manifold (SM), $\mathcal{L}\epsilon$, and the invariant manifold $\mathcal{M}_0$ perturbs into the adiabatic SSM (aSSM) denoted by $\mathcal{M}\epsilon$. This manifold $\mathcal{M}_\epsilon$ remains attracting and smoothly varies with the slow forcing. Its internal dynamics capture the essential, low-dimensional behavior of the full system, thereby serving as an accurate reduced-order model. Importantly, this reduced-order model can be learned directly from data, as we discuss in the Methods \ref{method_data_driven_aSSM} section.

In this adiabatic setting, we identify three well-separated time scales:
\begin{itemize}
    \item a \textit{fast time scale} $\tau_1$, associated with transient dynamics that lie outside the aSSM, $\mathcal{M}_\epsilon$;
    \item an \textit{intermediate time scale} $\tau_2$, governing the dynamics confined to the aSSM, $\mathcal{M}_\epsilon$;
    \item a \textit{slow time scale} $\tau_3$, determined by the rate of change of the external forcing, which also dictates the evolution on the slow manifold $\mathcal{L}_\epsilon$ to which $\mathcal{M}_\epsilon$ is attached.
\end{itemize}
These three time scales satisfy the hierarchy $\tau_1 < \tau_2 < \tau_3$, as schematically illustrated in Fig. \ref{adiabatic_cases}. Experimental observations of HASEL actuators and HASEL-based artificial muscles confirm this time scale hierarchy: the system dynamics exhibit rapid convergence to the aSSM, $\mathcal{M}_\epsilon$, followed by a slower alignment with the slow manifold, $\mathcal{L}_\epsilon$. Physically, this implies that internal transients decay much more rapidly than the evolution along the desired trajectory, i.e., $\tau_1 < \tau_2 \ll \tau_3$, as seen in Fig. \ref{adiabatic_cases}c. As a result, modeling the dynamics solely on the slow manifold provides an accurate and computationally efficient representation of the dominant dynamics of the system (see Method \ref{method_SM}), without resolving fast transients. Hence, the SM reduction offers a pragmatic solution for control, as the faster aSSM dynamics contribute only marginally to the performance.

\begin{figure*}[!htbp]
    \centering
    \includegraphics[scale=1]{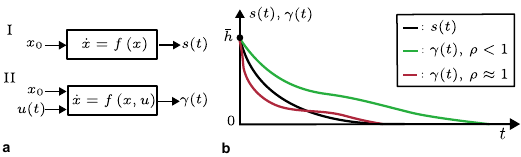}
    \caption{(\textbf{a}) Autonomous system ($\text{I}$) and externally forced system ($\text{II}$). (\textbf{b}) Comparison of trajectories decaying to zero from the same initial point $\bar{h}$ through autonomous transient decay and under slower ($\rho < 1$) and faster ($\rho \approx 1$) control in a regulation task.}
    \label{rho_visual}
\end{figure*}

We now introduce the non-dimensional metric $\rho$ to quantify the slowness of a control input signal $\boldsymbol{u}(t)$ (associated with the time scale $\tau_3$) relative to the autonomous decay $\boldsymbol{x}(t)$ (associated with $\tau_1$ and $\tau_2$), motivated by similar quantities defined in \citet{haller_2024} and \citet{kaundinya_2025}. We assume that a scalar observable $s(t)$ characterizes the autonomous system response and $\gamma(t)$ the response under external forcing $u(t)$. For a control application, $\gamma(t)$ would be the desired path. The observable $s(t)$ is shifted, if necessary, to vanish at the $x = 0$ origin. We further normalize $s(t)$ as
\begin{equation}
	\tilde{s}(t) = \frac{s(t)}{\left(\int_0^\infty |s(t)|^2 dt\right)^{1/2}}
\end{equation}
and $\gamma(t)$ as
\begin{equation}
	\tilde{\gamma}(t) = \frac{\gamma(t)}{\left(\int_{t_i}^{t_f} |\gamma(t)|^2 dt\right)^{1/2}},
\end{equation}
where the forced signal is considered in the time interval $[t_i, t_f]$. We define the dimensionless slowness measure
\begin{equation}
\rho = \frac{\left(\int_{t_i}^{t_f} |\dot{\tilde{\gamma}}(t)|^2 dt\right)^{1/2}}{\left(\int_0^\infty |\dot{\tilde{s}}(t)|^2 dt\right)^{1/2}}.
\end{equation}
For instance, for a linear autonomous system with observed displacement $s(t) = s_0 e^{\lambda t}$ and $\lambda < 0$, $\rho$ simplifies to
\begin{equation}
\rho = \frac{\left(\int_{t_i}^{t_f} |\dot{\tilde{\gamma}}(t)|^2 dt\right)^{1/2}}{| \lambda |}.
\end{equation}
In that case, the condition $\rho \ll 1$ ensures that the forced dynamics evolve slowly relative to the intrinsic decay of the system, thereby justifying the reduction to the slow manifold (see Fig. \ref{adiabatic_cases}c). Figure \ref{rho_visual} intuitively illustrates the meaning of the metric $\rho$ in a regulation task, comparing slow and fast control actions to the system's autonomous transient decay.

In the following, we adopt the normalized mean trajectory error (NMTE) from \citet{cenedese2022} as a measure of the average deviation between the predicted and actual trajectories to assess the accuracy of the model. The NMTE is defined as
\begin{equation}
	\text{NMTE} = \frac{1}{\|\boldsymbol{x}\|} \frac{1}{N}\sum_{j = 1}^N \|\boldsymbol{x}_j - \boldsymbol{\hat{x}}_j \|,
\end{equation}
where $\boldsymbol{x}$ denotes the true data, $\boldsymbol{\hat{x}}$ the corresponding model prediction, and $N$ the total number of temporal observations.

\subsection{Illustration of data-driven aSSM and SM modeling on analytic models}\label{analytical_models}
We first illustrate the power of aSSM-reduced modeling on data generated by analytic HASEL actuator models.
\subsubsection{Simple phenomenological model of HASEL actuator}
We consider a single-degree-of-freedom nonlinear oscillator that captures the qualitative behavior of HASEL actuators under voltage-driven excitation. This oscillator model is of the nondimensional form 
\begin{equation}\label{SDOF_oscillator}
m \ddot{x} + c\left(u\right) \dot{x} + k x + \alpha x^3 = \gamma u^2(t),
\end{equation}
with $c\left(u\right) = \tilde{c} - \beta u^2$ and with $u(t)$ representing the externally applied voltage. The term proportional to $\beta$ introduces a voltage-dependent damping mechanism, allowing the input to modulate the system’s effective relaxation time. Although the relationship between damping and voltage may be nonlinear and non-monotonic in experimental systems, the present model is designed to capture the potential influence of the input on internal decay rates. The term proportional to $\gamma$ acts as a purely external forcing input. To reproduce the overdamped behavior observed in experimental HASEL actuators, we set the model parameters to $m = 0.022$, $k = 1$, $\tilde{c} = 0.3$, $\alpha = 0.7$, $\beta = 5 \cdot 10^{-3}$, and $\gamma = 0.5$.

\begin{figure*}[!htbp]
    \centering
    \includegraphics[scale=1]{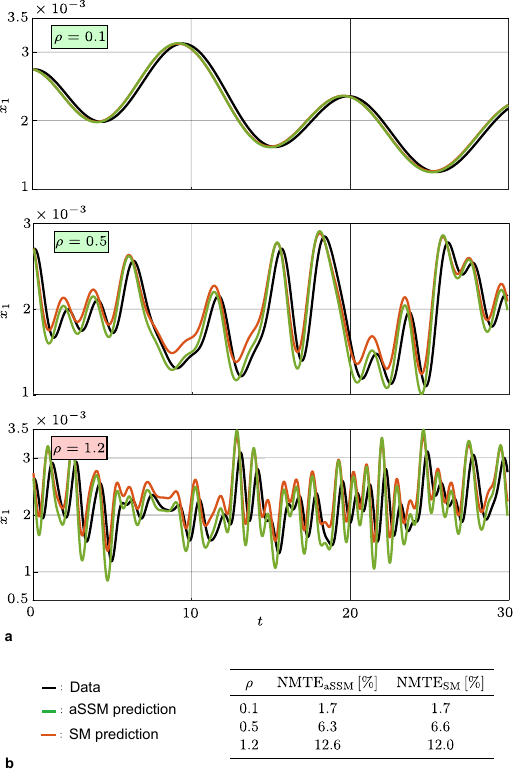}
    \caption{(\textbf{a}) Comparison between the aSSM-based and the SM-based reduced-order model of the SDOF of Eq. \eqref{SDOF_oscillator} under random external excitations with increasing speed ($\rho = 0.1$ (top), $\rho = 0.5$ (middle) and $\rho = 1.2$ (bottom)). As predicted by theory, the slow manifold approximation is accurate when $\rho \ll 1$, while in the last case even the leading-order aSSM-based approximation breaks down. (\textbf{b}) Legend and table reporting the NMTE values for the different cases.}
    \label{SDOF_model}
\end{figure*}

This example aims to illustrate how the rate of external forcing affects the predictive performance of reduced-order models (in this case one-dimensional) based on the aSSM and the SM approximation. We employ a leading-order approximation of the aSSM, as outlined in the Methods \ref{method_data_driven_aSSM}. The aSSM is constructed by collecting decaying trajectories at fixed voltage levels. For each fixed input, the system converges to a locally attracting SSM, from which the reduced dynamics and parametrization are learned. These local SSMs are then interpolated across voltage values to approximate the aSSM, $\mathcal{M}_\epsilon$ and its internal dynamics under slowly varying forcing. As a simplification of this procedure, the SM approach directly infers the input-state relationship from forced trajectories, without requiring decay experiments. This method assumes that internal dynamics decay sufficiently fast for the system state to remain close to a low-dimensional manifold defined instantaneously by the input.

Figure \ref{SDOF_model} compares the responses predicted by the aSSM and SM models under random voltage forcing. Although, in general, the SM is defined by a nonlinear function of the input, a linear approximation suffices in this simple example. We consider two cases with different forcing rates. For the fast forcing case ($\rho = 0.5$), the correction introduced by the leading-order aSSM significantly improves prediction accuracy. In contrast, as reported in Fig.~\ref{SDOF_model}b, when the forcing varies more slowly ($\rho = 0.1$), the system remains close to the SM, and the SM model yields sufficiently accurate predictions. Finally, for faster forcing ($\rho = 1.2$), both the SM- and aSSM-based predictions lose accuracy, as expected. Although such a high value of $\rho$ would be uncommon in practice, it illustrates the necessity of a time scale separation in the mathematical results underpinning aSSM reduction. We also emphasize that only the leading-order approximation of both the parametrization and the reduced dynamics of the  aSSM is used here. Higher-order approximations derived by \citet{haller_2024} would further enhance accuracy but would also require more data in order to avoid an overfit.


\subsubsection{Analytical model of HASEL actuator from the literature}\label{model_literature}

\begin{figure*}[!htbp]
    \centering
    \includegraphics[scale=1]{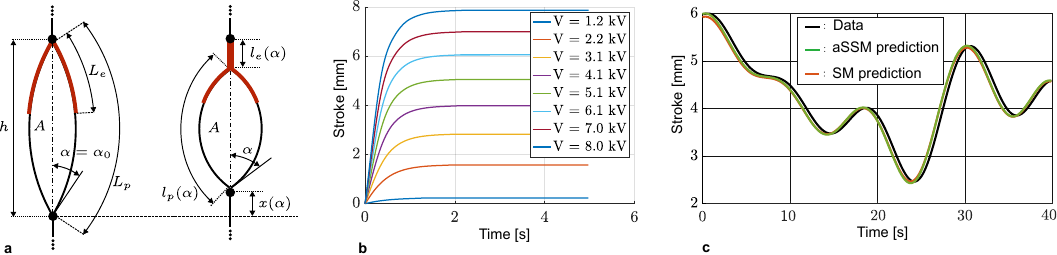}
    \caption{(\textbf{a}) Geometry of a single HASEL pouch actuator in the inactive state (left) and actuated state (right). Adapted from \citet{Kellaris_2019}. (\textbf{b}) Autonomous response of the analytic HASEL actuator model \eqref{governing_eq_analytic_HASEL} under different constant voltage levels. (\textbf{c}) Predictions of the data-driven aSSM-based and SM-based models for the actuator stroke in system \eqref{governing_eq_analytic_HASEL} under random excitation.}
    \label{analytic_geometry_results}
\end{figure*}

We now extend the simple model \eqref{SDOF_oscillator} of the HASEL actuator by incorporating the pouch geometry shown in Fig. \ref{analytic_geometry_results}a. This extension couples the deflection angle $\alpha$ of the pouch with the generated stroke $x$, and explicitly models the electrostatic forces acting on the electrodes. The pouches are treated as two-dimensional, assuming negligible influence from the out-of-plane depth. The actuator geometry follows the configuration presented in \citet{Kellaris_2019}. While we adopt a Lagrangian approach to derive the governing equations, these are equivalent to those obtained via the Port-Hamiltonian formalism  \citep{vanderSchaft_2014,Lohmayer_2021} employed by \citet{Yeh_2022}. Further details on the derivation are provided in the Supplementary Material.

The actuator consists of $N$ pouches connected in series. Each pouch is formed by bonding two rectangular dielectric films, each of length $L_p$, width $w$, and thickness $t$, along their edges. Within the experimental parameter range of interest, the dielectric films can be modeled as inextensible membranes without bending stiffness. Under these assumptions, the liquid-filled pouch adopts a cross-sectional shape composed of two intersecting circular segments with central angle $2\alpha_0$ and arc length equal to the undeformed pouch length $L_p$. The volume of fluid determines the cross-sectional area $A$, which is related to the angle $\alpha$ through
\begin{equation}
A = \frac{1}{2}L_p^2 \left(\frac{\alpha - \sin \alpha \cos \alpha}{\alpha^2}   \right).
\end{equation}
The initial actuator length is given by
\begin{equation}
h = L_p \left( \frac{\sin \alpha_0}{\alpha_0} \right).
\end{equation}
Electrodes of length $L_e$ and width $w$ are attached to the top surface of the actuator on both sides. When the voltage $u$ is applied, the electrodes zip together over a length $l_e$, starting from the edge of the pouch. If we neglect the thickness contribution from the dielectric fluid, the electrode separation becomes $2t$. Since the membrane is inextensible and the dielectric fluid is incompressible, the geometric parameters $L_p$, $w$, and $A$ remain constant throughout the actuation, leading to a vertical contraction of the actuator. During this process, the fluid-filled region retains the cylindrical segment shape, parametrized by the central angle $\alpha$.

It is possible to relate the arc length of the unzipped portion of the pouch $l_p$, the zipped length $l_e$, the stroke of the actuator $x$ and the capacitance of the zipped region $C$ to the angular deflection  $\alpha$ (see Supplementary Material). The model also incorporates a threshold voltage, below which the electrostatic force is insufficient to induce any zipping. Furthermore, it is assumed that the electrical energy arises solely from the region where the electrodes are zipped together, since the electric field in the unzipped region decays rapidly.

In terms of the stroke $x_i$ and the charge $Q_i$ on the zipped region of the electrodes of the $i$-th pouch, the governing equations read
\begin{equation}
\label{governing_eq_analytic_HASEL}
\begin{array}{l}
\displaystyle
m \ddot{x}_i + c \dot{x}_i + k x_i = mg - \frac{1}{2}\frac{Q_i^2 C_{x_i}}{C(x_i)^2}, \qquad i = 1,\dots, N \\[1ex] \displaystyle
\dot{Q}_i = -\frac{1}{RC(x_i)}Q_i + \frac{1}{R}u,
\end{array}
\end{equation}
where $g$ is the gravitational acceleration constant and $C_{x_i} = \frac{dC}{dx_i}$ describes the dependence of the capacitance on the actuator's stroke. The parameters of the model are tuned to mimic the behavior of the experimental HASEL actuator. 

\subsection{Data-driven modeling of experimental artificial muscles}\label{data_driven_modeling}

\begin{figure*}[!htbp]
    \centering
    \includegraphics[scale=1]{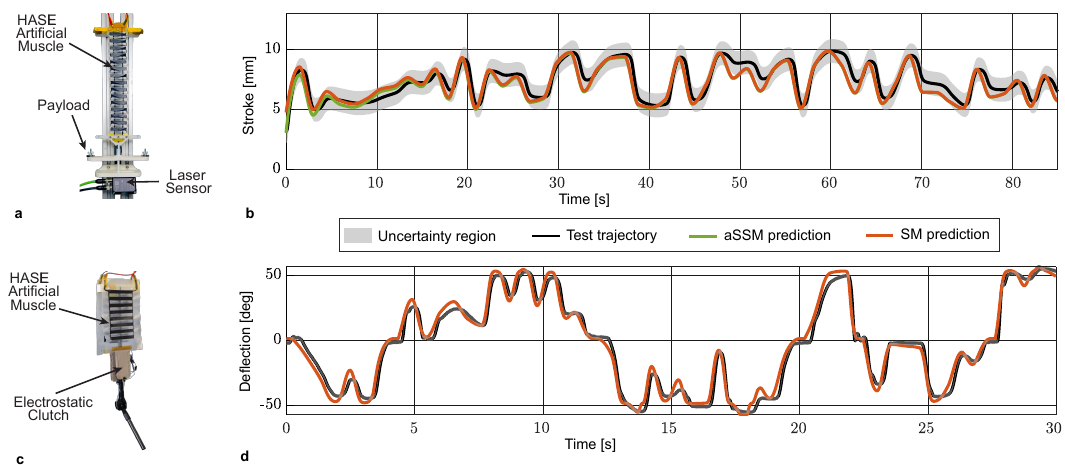}
    \vspace{1pt}
    \caption{\textbf{Experimental HASEL artificial muscles and data-driven slow manifold predictions.}
(\textbf{a}) Single HASEL actuator with payload and laser displacement sensor.
(\textbf{b}) Predictions of the actuator stroke under random forcing: the aSSM-based prediction (green) quickly converges to the SM-based prediction (orange), both closely matching the test trajectory (black).
(\textbf{c}) Musculoskeletal joint with an antagonistic pair of HASEL actuators and electrostatic clutches.
(\textbf{d}) Slow-manifold prediction of the joint deflection under random excitation, demonstrating accurate tracking of the test trajectory.}
    \label{HASEL_experiments}
\end{figure*}

We now consider the experimental HASEL actuator shown in Fig. \ref{HASEL_experiments}a, which consists of fifteen pouches. Experimental results under random voltage forcing reveal that the actuator's intrinsic dynamics are significantly faster than the typical rates of change of the inputs necessary for producing typical HASEL contraction patterns arising in practice. This pronounced time-scale separation justifies the use of a one-dimensional reduced-order model based on a slow manifold, given that the input $u(t)$ is one-dimensional.

To build such a model, we use the stroke of the actuator as a single observable and train an SM-based reduced model directly from forced data. The SM is approximated by a polynomial function trained on seven randomly forced trajectories. Notably, a linear approximation
\begin{equation}
x_\epsilon\left(u\right) = 1.12 \cdot 10^{-7} + 1.56 \cdot 10^{-3} u,
\end{equation}
proves remarkably effective for this system. For comparison, we also construct an aSSM-based model using a cubic-order approximation for both the manifold parametrization and the reduced dynamics. As shown in Fig. \ref{HASEL_experiments}b, the prediction from the aSSM model rapidly converges to that of the SM model, offering no significant improvement. This behavior remains consistent across voltage inputs with different slowness measure $\rho$. The example shown in Fig. \ref{HASEL_experiments}b features the dimensionless rate $\rho = 0.28$, yielding normalized mean trajectory errors of $\mathrm{NMTE}_\mathrm{aSSM} = 7.90\%$ and $\mathrm{NMTE}_\mathrm{SM} = 8.16\%$. As shown in the figure, both errors fall within the range of experimental uncertainty, estimated from repeated trials using the same voltage inputs.

The construction of the aSSM model, as discussed in \citet{kaundinya_2025}, requires generating a set of initial conditions as perturbations from steady states corresponding to different constant voltages. For a HASEL actuator, these new initial conditions are produced by applying sudden voltage jumps relative to the constant input. Such high-amplitude jumps, however, often trigger hysteretic actuator behavior, that would otherwise not arise under typical operating conditions. This mismatch limits the applicability of the aSSM approach based on experimental data for this class of systems. 

Each pouch is described by Eq. \eqref{governing_eq_analytic_HASEL}, which corresponds to three first-order ordinary differential equations, so the phase-space dimension of the entire actuator is $n = 3N$. The dimension of the external forcing (voltage) is $f = 1$. Given that the actuator exhibits overdamped behavior (see Fig. \ref{analytic_geometry_results}b) and the dynamics of the charge on the electrodes evolve much faster, the reduced-order model is effectively one-dimensional. The single observable used for data-driven modeling is the total stroke of the actuator. Figure \ref{analytic_geometry_results}c presents a comparison between the predictions of the aSSM-based and SM-based reduced-order models under random forcing with $\rho = 0.13$. The models yield normalized mean trajectory errors of $NMTE_{aSSM} = 1.7 \%$ and $NMTE_{SM} = 1.8 \%$, respectively. These results confirm the accuracy of data-driven aSSM and SM models for system \eqref{governing_eq_analytic_HASEL}.

In contrast, the SM-based reduced modeling relies solely on forced response data, which can be collected under normal operating conditions without introducing unwanted hysteresis. This practical advantage, combined with the comparable prediction accuracy of the SM model, makes the slow manifold approach more suitable for modeling HASEL actuators. Further details on the actuator’s hysteresis, as well as on data collection and processing, can be found in the Supplementary Material.

We now extend our analysis from a single actuator to a complete artificial muscle system. Specifically, we investigate an antagonistic configuration composed of two HASEL actuators arranged on opposite sides of a rigid structure, each paired with an electrostatic clutch. This setup, inspired by the design presented by \citet{Kazemipour_2024} and shown in Fig. \ref{HASEL_experiments}c, enables an extended range of motion and allows the artificial muscle to emulate the antagonistic actuation of biological limbs.

In this configuration, the system features four control inputs: two for activating the HASELs and two for engaging the electrostatic clutches. In the present work, we assume that the HASEL and clutch on a given side are activated simultaneously, and that only one side of the muscle is active at a time. This assumption reduces the input space to a single scalar input: the sign of this input determines which side of the muscle is engaged, thus preserving the antagonistic actuation behavior in a simplified representation.

Unlike isolated actuators, designing decaying trajectories for artificial muscles is significantly more challenging due to the complexity of the multi-input configuration and mechanical coupling. This makes the use of aSSMs less practical, as their construction requires decaying data whose generation induces undesirable hysteretic behavior that is atypical under normal operating conditions. In contrast, the SM-based approach demonstrates its practicality by requiring only forced trajectories for training. In this case, we employ a polynomial approximation of the slow manifold up to 7th order, trained using just three forced-response trajectories. As shown in Fig. \ref{HASEL_experiments}d, the SM-based model accurately predicts the muscle deflection, yielding a normalized mean trajectory error of $\mathrm{NMTE}_\mathrm{SM} = 8.13\%$.\\

Having validated the predictive accuracy of the data-driven slow manifold models, we now demonstrate their utility for real-time control by implementing model-based feedforward on the antagonistic joint and comparing performance against feedback-only and feedforward-only baselines.

\subsection{SSM-enabled closed-loop control of an antagonistic joint}\label{control_section}
\begin{figure*}[!htbp]
    \centering
    \includegraphics[scale=1]{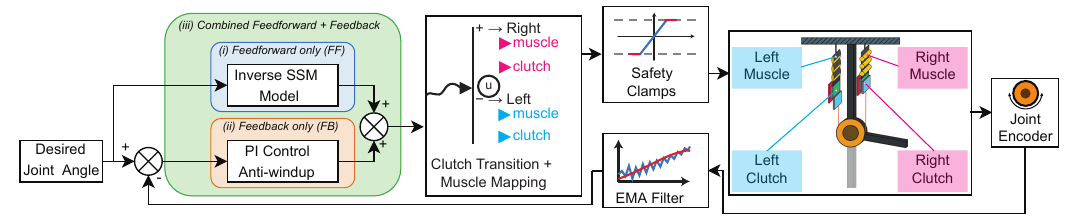}
    \caption{Closed-loop control architecture (FF, PI, FF+PI). (i) Feedforward: the desired joint angle is mapped to a voltage command by the inverse SSM model \(g^{-1}\). (ii) Feedback: a proportional--integral (PI) controller with anti-windup acts on the exponential moving average (EMA)-filtered joint-encoder measurement. (iii) Combined controller: feedforward and feedback are summed and routed through the clutch-transition and muscle-mapping stage; safety clamps then enforce saturation and slew-rate limits before driving the antagonistic left/right muscles and their clutches.}
    \label{fig:control_arch}
\end{figure*}

\begin{figure*}[!htbp]
    \centering
    \includegraphics[scale=1]{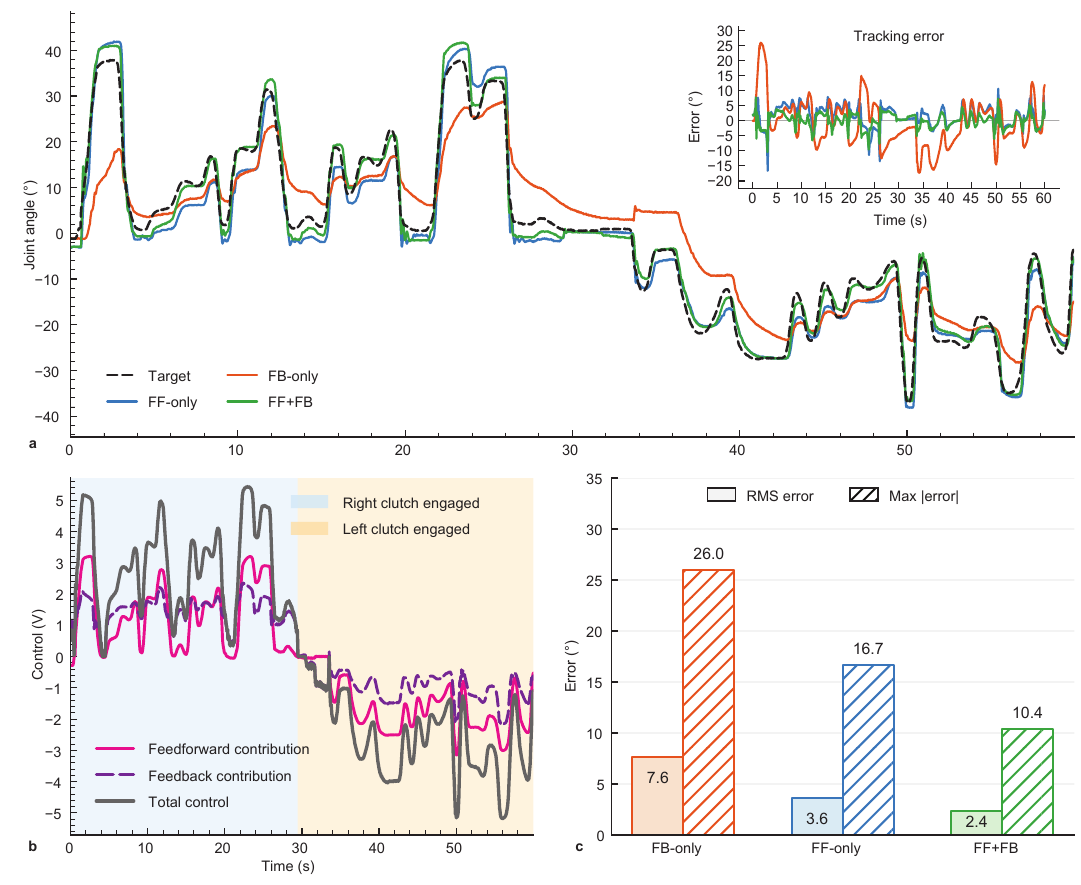}
    \caption{\textbf{Combined model-based feedforward and feedback control achieves superior tracking performance.} (\textbf{a}) Target joint angle (dashed) and measured angle for three controllers: PI-only, inverse-SSM feedforward-only, and combined inverse-SSM+PI. (\textbf{b}) Control signal for the combined controller, decomposed into feedforward and feedback contributions; shaded bands indicate engagement of the right and left clutches. (\textbf{c}) Performance summary showing root-mean-square (RMS) error and maximum absolute error \(|e|_{\max}\), demonstrating two-thirds error reduction relative to feedback-only control.}
    \label{fig:tracking_cases}
\end{figure*}

In order to test our data-driven SM-model in a closed-loop control setting, we observe a time-varying joint angle \(\theta_d(t)\) using the antagonistic HASEL--clutch setup of Fig.~\ref{HASEL_experiments}c, with measurement \(\theta(t)\). We evaluate three control strategies: proportional–integral (PI) feedback, an inverse-SSM feedforward term added to PI, and inverse-SSM feedforward alone. The PI law is \(u(t)=k_p e(t)+k_i \int_0^t e(\tau)\,\mathrm{d}\tau\) with \(e(t)=\theta_d(t)-\theta(t)\). The inverse-SSM term is obtained from the identified slow-manifold map \(\theta\approx g(u)\) (Section~\ref{data_driven_modeling}) via \(g^{-1}(\theta_d(t))\).

The gains \((k_p,k_i)\) are tuned separately for the combined and feedback-only controllers on a short calibration trajectory and then held fixed. All controllers are run on a distinct, randomized reference under identical safety constraints and clutch scheduling. Performance is quantified by the root-mean-square (RMS) tracking error and the maximum absolute error, \(\max_t |e(t)|\).

On the unseen trajectory, the inverse-SSM+PI controller achieved the lowest error (RMS 2.38$^\circ$; max 10.41$^\circ$). Feedforward-only attained RMS 3.63$^\circ$ and max 16.67$^\circ$, while feedback-only produced RMS 7.63$^\circ$ and max 26.00$^\circ$. The control effort, quantified directly from the commanded voltages to the two muscles (\(|u_\text{right}|+|u_\text{left}|\) per sample), was comparable across conditions in RMS and mean magnitude (mean \(\approx\) 2.68--2.98\,V; RMS \(\approx\) 3.01--3.15\,V), with feedback-only exhibiting the largest mean magnitude (2.98\,V). Figure~\ref{fig:tracking_cases}a overlays the target and the measured angles for all conditions; Fig.~\ref{fig:tracking_cases}b decomposes the combined control into feedforward and feedback contributions with shaded clutch engagement; Fig.~\ref{fig:tracking_cases}c summarizes the RMS and maximum errors.

When used alone, the inverse-SSM control strategy tracks the slow trend but cannot reject disturbances, leading to higher peaks in \(|e|\). The PI-only strategy corrects errors but requires larger voltages and exhibits slower transients when the antagonistic side switches. In the combined control strategy, the commanded-voltage contributions confirm that feedforward dominates effort: the feedforward term has mean magnitude 1.47\,V (RMS 1.77\,V) versus 1.21\,V (RMS 1.32\,V) for the feedback term, i.e., a \(\approx\)1.22\,\(\times\) larger mean and \(\approx\)1.34\,\(\times\) larger RMS contribution. 

Taken together, these observations indicate that the SSM-based feedforward supplies the dominant actuation burden, while the PI term effectively rejects disturbances and residual model error.

\section{Discussion}\label{discussion}

We presented a data-driven reduced-order modeling methodology for artificial muscles and demonstrated its performance on closed-loop control of an antagonistic HASEL–clutch joint. The central observation across our simulations and experiments is a pronounced time‑scale separation required for aSSM-based model reduction: trajectories collapse rapidly onto an attracting, low‑dimensional set before evolving along a slowly varying response. In this regime, a slow‑manifold (SM) model learned directly from forced data suffices for prediction and control, while a full aSSM reduction provides additional accuracy. The latter, however, offers modest additional benefit in practice because it requires fast decay experiments that can excite hysteresis. The advantage of using the reduction on the aSSM instead of the SM would become more evident under milder time-scale separation, for instance when the forcing varies more rapidly or in soft robots exhibiting slower internal dynamics. 

Importantly, the data-driven approach modeling itself is our main contribution: the SM/aSSM workflow developed here is sample‑efficient (trained from a handful of forced trajectories), predictive (single‑digit NMTE on both single‑actuator and antagonistic‑muscle tests), and computationally light. It also yields an explicit inverse map \(g^{-1}\) with microsecond‑level evaluation, enabling real‑time feedforward on embedded hardware.

A key outcome of this work is demonstrating that data-driven reduced-order models are not merely accurate predictors but also effective control tools for antagonistic artificial muscle systems. Antagonistic actuation, inspired by biological muscle pairs, offers variable impedance and bidirectional motion, yet controlling such assemblies has remained challenging due to switching dynamics, multi-physics coupling, and the need for real-time computation. By successfully implementing model-based feedforward control on an antagonistic HASEL-clutch joint and achieving substantial error reduction over feedback-only baselines, we establish that the learned slow-manifold inverse provides a practical path to precise soft-robot control without prohibitive computational overhead. This bridges a critical gap between theoretical model reduction and deployable robotic systems, showing that bio-inspired antagonistic designs can be controlled accurately using computationally lean, experiment-derived models.

Our control experiments have substantiated these modeling choices. Using the learned inverse map \(g^{-1}\) as feedforward and a lightweight PI feedback term, the combined controller achieves the lowest tracking error on an unseen randomized trajectory (Fig.~\ref{fig:tracking_cases}). Relative to feedforward‑only and PI‑only, the RMS error decreases by roughly one‑third and two‑thirds, respectively, while the maximum error improves accordingly. Notably, the commanded voltage magnitudes delivered to the muscles are of similar scale across all three conditions, indicating that accuracy gains are not obtained by substantially increasing actuation. Within the combined controller, the decomposition of the command shows that the feedforward contribution accounts for the larger share of the effort (mean magnitude \(\approx\)1.47\,V vs. \(\approx\)1.21\,V for feedback), with a comparable dominance in RMS. Thus, the model supplies most of the actuation burden, and the feedback primarily rejects disturbances and residual modeling error.

These findings suggest a practical recipe for soft‑robotic muscle control under slowly varying motion goals: (i) learn a one‑dimensional slow‑manifold map from operational data; (ii) use its inverse as a fast, robust feedforward term; and (iii) add a small PI loop with anti‑windup for robustness. The resulting controller is simple to tune, executes at kHz rates, and remains compatible with antagonistic actuation and clutch scheduling. Beyond the present system, the same modeling pipeline applies to other soft‑robotic actuators that display clear time‑scale separation, providing interpretable, portable models with direct control utility.

The demonstrated approach has implications beyond HASEL actuators. Other soft robotic systems that exhibit pronounced time-scale separation, such as pneumatic artificial muscles (PAMs), certain dielectric elastomer actuators (DEAs), or shape-memory alloy (SMA) actuators with fast thermal transients, may benefit from this sample-efficient modeling framework. The ability to learn accurate models from a handful of training trajectories makes this practical for field deployment where extensive data collection is infeasible, opening pathways for adaptive prosthetics, minimally invasive surgical robotics, and collaborative manufacturing where soft, precise actuation is essential. Moreover, the computationally lean nature of the inverse slow-manifold map enables deployment on resource-constrained embedded controllers typical of wearable and portable robotic systems.

Our approach is not without limitations. The SM approximation presumes slow external variation (small slowness rate \(\rho\)), and its fidelity degrades under abrupt voltage changes that elicit hysteresis. Our antagonistic realization assumes mutually exclusive clutch engagement; richer behaviors (e.g., brief overlap or bilateral actuation) will require higher‑dimensional slow manifolds. Finally, while the present inverse map is static, mild rate‑dependence could further improve accuracy under faster references.

Based on these limitations, future work will target (i) asynchronous HASEL/clutch activation and simultaneous bilateral actuation via higher‑dimensional aSSMs/SMs; (ii) hysteresis‑aware or rate‑dependent corrections for faster inputs; (iii) online adaptation of \(g\) to compensate drift; and (iv) principled optimization of clutch timing to balance performance and safety.

\section{Methods}\label{methods}

\subsection{Data-driven modeling via adiabatic spectral submanifolds}\label{method_data_driven_aSSM}
For $\epsilon = 0$, we assume that the system  \eqref{nonautonomous_DS} has a $\boldsymbol{u}$-dependent, uniformly bounded and uniformly hyperbolic fixed point $\boldsymbol{x}_0(\boldsymbol{u})$, as explained in \citet{haller_2024}. In this setting, the system has an $f$-dimensional, normally attracting invariant manifold 
\begin{equation}\label{critical_manifold}
\mathcal{L}_0 = \left\{\left(\boldsymbol{x}, \boldsymbol{u} \right) \in \mathbb{R}^n \times \mathbb{R}^f :   \boldsymbol{x} = \boldsymbol{x}_0\left( \boldsymbol{u}\right) \right\},
\end{equation}
which is the collection of the $\boldsymbol{x}_0\left(\boldsymbol{u}\right)$ fixed points and is called critical manifold in the terminology of the geometric singular perturbation theory \citep{fenichel_1979}. The manifold $\mathcal{L}_0$ then satisfies
\begin{equation}
\boldsymbol{A}\boldsymbol{x}_0\left(\boldsymbol{u}\right) + \boldsymbol{f}_0\left( \boldsymbol{x}_0\left(\boldsymbol{u}\right) \right) + \boldsymbol{f}_1\left(\boldsymbol{x}_0\left(\boldsymbol{u}\right), \boldsymbol{u} \right) = \boldsymbol{0}.
\end{equation}

We now perform a change of variables
\begin{equation}
\boldsymbol{\xi} = \boldsymbol{x} - \boldsymbol{x}_0\left( \boldsymbol{u}\right),
\end{equation}
which flattens out the critical manifold. We can then rewrite eq. \eqref{nonautonomous_DS}
\begin{equation}\label{reduced_dynamics_adiabatic}
\begin{cases}
\begin{aligned}
\dot{\boldsymbol{\xi}} &= \boldsymbol{A} \left( \boldsymbol{x}_0\left(\boldsymbol{u}\right)+\boldsymbol{\xi}\right) 
+ \boldsymbol{f}_0\left(\boldsymbol{x}_0\left(\boldsymbol{u}\right)+\boldsymbol{\xi}\right) \\
&\quad + \boldsymbol{f}_1\left(\boldsymbol{x}_0\left(\boldsymbol{u}\right)+\boldsymbol{\xi}, \boldsymbol{u}\right)
- \boldsymbol{D}_{\boldsymbol{u}}\boldsymbol{x}_0\left(\boldsymbol{u}\right)\dot{\boldsymbol{u}}, \\
\dot{\boldsymbol{u}} &= \epsilon \boldsymbol{w}\left( \boldsymbol{x}_0\left(\boldsymbol{u}\right)+\boldsymbol{\xi}, \epsilon t \right).
\end{aligned}
\end{cases}
\end{equation}

We now expand in Taylor series the right-end-side of the first equation in \eqref{reduced_dynamics_adiabatic} to obtain 
\begin{equation}\label{aSSM_equations_shifted}
\begin{split}
\boldsymbol{f}_0\left(\boldsymbol{x}_0\left(\boldsymbol{u}\right)+\boldsymbol{\xi}\right) 
&= \boldsymbol{f}_0\left(\boldsymbol{x}_0\left(\boldsymbol{u}\right)\right) 
+ \boldsymbol{D}_{\boldsymbol{x}}\boldsymbol{f}_0\left(\boldsymbol{x}_0\left(\boldsymbol{u}\right)\right)\boldsymbol{\xi} \\
&\quad + \frac{1}{2}\boldsymbol{\xi}^\mathrm{T} \boldsymbol{D}^2_{\boldsymbol{x}}\boldsymbol{f}_0\left(\boldsymbol{x}_0\left(\boldsymbol{u}\right)\right)\boldsymbol{\xi} 
+ \mathcal{O}\left(3\right),
\end{split}
\end{equation}
and 
\begin{equation}
\begin{split}
\boldsymbol{f}_1\left(\boldsymbol{x}_0\left(\boldsymbol{u}\right)+\boldsymbol{\xi}, \boldsymbol{u}\right) 
&= \boldsymbol{f}_1\left(\boldsymbol{x}_0\left(\boldsymbol{u}\right), \boldsymbol{u}\right) 
+ \boldsymbol{D}_{\boldsymbol{x}}\boldsymbol{f}_0\left(\boldsymbol{x}_0\left(\boldsymbol{u}\right), \boldsymbol{u}\right)\boldsymbol{\xi} \\
&\quad + \frac{1}{2}\boldsymbol{\xi}^\mathrm{T} \boldsymbol{D}^2_{\boldsymbol{x}}\boldsymbol{f}_1\left(\boldsymbol{x}_0\left(\boldsymbol{u}\right),\boldsymbol{u}\right)\boldsymbol{\xi} 
+ \mathcal{O}\left(3\right).
\end{split}
\end{equation}
Equation \eqref{reduced_dynamics_adiabatic} thus becomes
\begin{equation}
\begin{cases}
\begin{aligned}
\dot{\boldsymbol{\xi}} &= \hat{\boldsymbol{A}}\left(\boldsymbol{u}\right) \boldsymbol{\xi}  - \boldsymbol{D}_{\boldsymbol{u}}\boldsymbol{x}_0\left( \boldsymbol{u} \right) \dot{\boldsymbol{u}} + \frac{1}{2} \boldsymbol{\xi}^\mathrm{T} \left\{ \boldsymbol{D}_{\boldsymbol{x}}\boldsymbol{f}_0\left( \boldsymbol{x}_0\left( \boldsymbol{u} \right) \right) \right. \\
&\quad\quad \left. + \boldsymbol{D}^2_{\boldsymbol{x}}\boldsymbol{f}_1\left( \boldsymbol{x}_0\left( \boldsymbol{u} \right), \boldsymbol{u} \right) \right\} \boldsymbol{\xi} + \mathcal{O}\left(3\right), \\
\dot{\boldsymbol{u}} &= \epsilon \boldsymbol{w}\left( \boldsymbol{x}_0\left( \boldsymbol{u} \right) + \boldsymbol{\xi}, \epsilon t \right),
\end{aligned}
\end{cases}
\end{equation}
with the linear part
\begin{equation*}
\hat{\boldsymbol{A}}\left( \boldsymbol{u}\right)= \boldsymbol{A} + \boldsymbol{D}_{\boldsymbol{x}}\boldsymbol{f}_0\left( \boldsymbol{x}_0\left( \boldsymbol{u} \right) \right) 
+ \boldsymbol{D}_{\boldsymbol{x}}\boldsymbol{f}_1\left( \boldsymbol{x}_0\left( \boldsymbol{u} \right), \boldsymbol{u} \right) \in \mathbb{R}^{n \times n}.
\end{equation*}
We will denote by $\boldsymbol{P}\left(\boldsymbol{u}\right) \in \mathbb{R}^{n \times n}$ the column matrix of the right eigenvectors of $\hat{\boldsymbol{A}}\left( \boldsymbol{u}\right)$ and $\boldsymbol{Q}_{\boldsymbol{u}}\left(\boldsymbol{u}\right) \in \mathbb{R}^{d \times n}$ the row matrix of the first $d$ left eigenvectors. 

We now apply a modal coordinate transformation 
\begin{equation}
\boldsymbol{\xi} = \boldsymbol{P}\left( \boldsymbol{u} \right) \begin{pmatrix}
\boldsymbol{\eta} \\ \boldsymbol{v}
\end{pmatrix}, \quad \boldsymbol{\eta} \in \mathbb{R}^d,\quad \boldsymbol{v} \in \mathbb{R}^{n - d}.
\end{equation}
The leading-order approximation of the reduced dynamics on the adiabatic SSM then becomes  
\begin{equation}\label{aSSM_reduced_dynamics}
\begin{aligned}
\dot{\boldsymbol{\eta}} &= 
\boldsymbol{Q}_{\boldsymbol{u}} \boldsymbol{\hat{A}}\left(\boldsymbol{u}\right) 
\boldsymbol{P}\left( \boldsymbol{u} \right)
\begin{pmatrix}
\boldsymbol{\eta} \\
\boldsymbol{h}_0\left(\boldsymbol{\eta}, \boldsymbol{x}_0\left(\boldsymbol{u}\right)\right)
\end{pmatrix} \\&\quad
+ \frac{1}{2}
\begin{pmatrix}
\boldsymbol{\eta} \\
\boldsymbol{h}_0\left(\boldsymbol{\eta}, \boldsymbol{x}_0\left(\boldsymbol{u}\right)\right)
\end{pmatrix}^\mathrm{T} 
\boldsymbol{P}\left( \boldsymbol{u} \right)^\mathrm{T}
\Biggl(
\boldsymbol{D}_{\boldsymbol{x}}\boldsymbol{f}_0\left( \boldsymbol{x}_0\left( \boldsymbol{u} \right) \right)\\ &\quad
+ \boldsymbol{D}^2_{\boldsymbol{x}}\boldsymbol{f}_1\left( \boldsymbol{x}_0\left( \boldsymbol{u} \right), \boldsymbol{u} \right)
\Biggr) \boldsymbol{P}\left( \boldsymbol{u} \right)
\begin{pmatrix}
\boldsymbol{\eta} \\
\boldsymbol{h}_0\left(\boldsymbol{\eta}, \boldsymbol{x}_0\left(\boldsymbol{u}\right)\right)
\end{pmatrix}
+ \mathcal{O}\left(3\right),
\end{aligned}
\end{equation}
where the relationship $\boldsymbol{v} = \boldsymbol{h}_0\left(\boldsymbol{\eta},  \boldsymbol{x}_0(\boldsymbol{u})\right)$ is the leading-order approximation in $\epsilon$ of the aSSM parametrization. For a discussion on higher-order corrections, see \citet{haller_2024}.

Given a scalar observable $s(t) \in \mathbb{R}^N$, such as the stroke of the HASEL actuator or the deflection angle of the musculoskeletal joint, we apply \textit{delay embedding} to ensure that the embedding dimension is sufficient to reconstruct the manifold in the observable space \citep{takens_1981}. The following procedure enables us to recover the leading-order approximation of the aSSM and its associated reduced dynamics for a scalar external input $u \in \mathbb{R}$:

\begin{enumerate}
\item Given the observed data $s \in \mathbb{R}^{N}$ at the frozen time limit $u = \bar{u} \in \mathbb{R}$, assess the steady solution $s_0(\bar{u})$.

\item Delay embed the data and construct $\boldsymbol{x}(t)$ and $\boldsymbol{x}_0(\bar{u})$, with $\boldsymbol{x}, \boldsymbol{x}_0(\bar{u}) \in \mathbb{R}^{p \times N}$.

\item Define the shift  $\boldsymbol{\xi} = \boldsymbol{x} - \boldsymbol{x}_0(\bar{u})$ with respect to the fixed point $\boldsymbol{x}_0(\bar{u})$.

\item Approximate the tangent space $\boldsymbol{V} \in \mathbb{R}^{p \times d}$ of $\boldsymbol{\xi}(t)$ via singular value decomposition (SVD), where $d$ is the dimension of the reduced-order model. In our case, $d = 1$.

\item Define the reduced coordinates $\boldsymbol{\eta} = \boldsymbol{V}^\mathrm{T} \boldsymbol{x}, \quad \boldsymbol{\eta} \in \mathbb{R}^{d \times N}$. 

\item Compute the SSM parametrization $\boldsymbol{\xi} = \boldsymbol{h}_0\left(\boldsymbol{\eta}; \bar{u} \right) = \left[ \boldsymbol{H}_1(\bar{u}), \dots \boldsymbol{H}_m(\bar{u}) \right] \boldsymbol{\eta}^{1:m}$, which is a multivariate polynomial and $\boldsymbol{H}_i \in \mathbb{R}^{p \times d_i}$ and it is the zeroth-order approximation in $\epsilon$ of the aSSM at $u = \bar{u}$.

\item Compute the reduced dynamics $\boldsymbol{\dot{\eta}} = \boldsymbol{r}\left(\boldsymbol{\eta}; \bar{u}\right) = \left[\boldsymbol{R}_1(\bar{u}),\dots,\boldsymbol{R}_r(\bar{u})  \right] \boldsymbol{\eta}^{1:r}$, which is again a multivariate polynomial, where $\boldsymbol{R}_i \in \mathbb{R}^{d \times d_i}$ and it is an approximation of Eq. \eqref{aSSM_reduced_dynamics} up to the $r^{th}$ order in the reduced coordinate $\boldsymbol{\eta}$.

\item Compute the SSM parametrization as
\begin{equation}
\boldsymbol{\xi} = \boldsymbol{h}_0(\boldsymbol{\eta}; \bar{u}) 
= \sum_{|\boldsymbol{k}|=1}^{m} \boldsymbol{H}_{\boldsymbol{k}}(\bar{u})\, \boldsymbol{\eta}^{\boldsymbol{k}},
\end{equation}
where $\boldsymbol{H}_{\boldsymbol{k}}(\bar{u}) \in \mathbb{R}^{p \times d_{\boldsymbol{k}}}$ are coefficient matrices, $\boldsymbol{k} = (k_1, \ldots, k_d) \in \mathbb{N}^d$ is a multi-index, 
$|\boldsymbol{k}| = k_1 + \cdots + k_d$ denotes its total degree, 
and $\boldsymbol{\eta}^{\boldsymbol{k}} = \eta_1^{k_1} \eta_2^{k_2} \cdots \eta_d^{k_d}$.
This expression represents a multivariate polynomial corresponding to the zeroth-order approximation in $\epsilon$ of the aSSM evaluated at $u = \bar{u}$.

\item Compute the reduced dynamics as
\begin{equation}
\dot{\boldsymbol{\eta}} = \boldsymbol{r}(\boldsymbol{\eta}; \bar{u}) 
= \sum_{|\boldsymbol{k}|=1}^{r} \boldsymbol{R}_{\boldsymbol{k}}(\bar{u})\, \boldsymbol{\eta}^{\boldsymbol{k}},
\end{equation}
where $\boldsymbol{R}_{\boldsymbol{k}}(\bar{u}) \in \mathbb{R}^{d \times d_{\boldsymbol{k}}}$ are coefficient matrices. 
This representation provides a multivariate polynomial approximation of Eq.~\eqref{aSSM_reduced_dynamics} up to order $r$ in the reduced coordinates $\boldsymbol{\eta}$.

\vspace{1ex}

\end{enumerate}
This procedure is repeated for each voltage value $u$ in the sampling. The resulting fixed points $\boldsymbol{x}_0(u)$ are then interpolated to reconstruct the geometry of the critical manifold, along with the associated parametrizations and reduced dynamics as functions of the voltage $u$. This interpolation assumes that the SSM dimension remains constant across the sampled range. While the tangent spaces of the SSM may vary with voltage, one may interpolate them or, to reduce computational cost, assume a constant tangent space. In the latter case, it is essential to ensure that this fixed plane remains sufficiently transverse to the fast subspace at all voltage values. If the fast subspace is strongly non-normal relative to the slow one, oblique projections may also be employed (see \citet{bettini_2025}).

\subsection{Data-driven modeling via slow manifold}\label{method_SM}

For $\epsilon > 0$ small enough, the critical manifold $\mathcal{L}_0$ as defined in Eq. \eqref{critical_manifold} perturbs into a unique, attracting, and $f$-dimensional slow manifold (SM) 

\begin{equation}\label{slow_manifold}
\mathcal{L}_0 = \left\{\left(\boldsymbol{x}, \boldsymbol{u} \right) \in \mathbb{R}^n \times \mathbb{R}^f :   \boldsymbol{x} = \boldsymbol{x}_\epsilon\left( \boldsymbol{u}\right) \right\},
\end{equation}
which is diffeomorphic to $\mathcal{L}_0$. For further details, we refer to \citet{haller_2024}.

As observed from experimental data and seen in Fig. \ref{adiabatic_cases}c, system trajectories rapidly converge not only to the aSSM, but also to the underlying SM under slowly varying external inputs. This observation enables the direct identification of the SM from forced response data, effectively treating it as a mapping from the input into the observable.
For a scalar observable $s(t) \in \mathbb{R}$ and a scalar input $u(t) \in \mathbb{R}$, the SM can be approximated as a polynomial function
\begin{equation}
x_\epsilon(u) \doteq s(u) = \sum_{i = 0}^n S_i u^i, \quad S_i \in \mathbb{R}.
\end{equation} 
Alternative approximations (radial basis functions, neural networks; see \citet{Powell_1992,Hagan_2014,Park_1991}) could also be used, but have turned out to be less efficient in our tests (see the Supplementary Material).

\subsection{Control design and implementation}\label{methods_control}

As noted in Section \ref{control_section}, we evaluate three controller realizations for the antagonistically actuated HASEL elbow: (i) feedforward only, (ii) feedback only, and (iii) combined feedforward+feedback. All three realizations respect identical safety constraints, apply the same actuation timing for the clutched antagonists, and run with a fixed sampling period \(T_s\) (1\,kHz in our implementation). The overall signal flow is summarized in the control block diagram (Fig.~\ref{fig:control_arch}), which comprises an inverse SSM-based feedforward path, a PI feedback path with anti-windup operating in the EMA-filtered measurement, a clutch transition and actuator mapping stage, and terminal safety clamps that enforce hardware limits before driving the antagonistic muscles and clutches. In the following, we summarize the control laws, discretization details, and the antagonistic actuation workflow.

\textit{Notation:} Let \(\theta_d(t)\) denote the desired joint angle and \(\theta(t)\) the measured angle. The control signal \(V(t)\) is the commanded high-voltage applied to the active muscle. Safety constraints include a voltage saturation \(V\in[\underline V,\overline V]\) and a slew-rate bound \(|V(t)-V(t-T_s)|\leq \Delta V_{\max}\). In discrete time, the filtered angle \(\tilde\theta_k\) is obtained by a first-order low-pass on the raw measurement to suppress quantization/noise.

\textit{Feedforward from the SSM-based model:} We identify a static, monotone map between voltage and quasi-static joint angle, \(\theta \approx g(V)\) (Section~\ref{data_driven_modeling}). The feedforward command inverts this relation at the desired angle trajectory:
\[
V_{\mathrm{ff}}(t) = g^{-1}\!\big(\theta_d(t)\big).
\]

\textit{Numerical inverse}: We precompute a dense one-dimensional interpolant over \((V,\,g(V))\) from the learned polynomial forward map and recover \(V_{\mathrm{ff}}\) by linear interpolation. This avoids per-sample root finding and yields microsecond-level runtime.
\textit{Range protection}: \(\theta_d\) is clipped to the calibrated image of \(g\) to avoid extrapolation; the resulting \(V_{\mathrm{ff}}\) is then clipped to \([\underline V,\overline V]\) and rate-limited by \(\Delta V_{\max}\).

In the feedforward-only benchmark, the applied command is
\[
V(t) = \operatorname{sat}_{\underline V,\overline V}\!\big(u(t)\big),
\qquad
u(t) = \mathrm{rl}_{\Delta V_{\max}}\{V_{\mathrm{ff}}(t)\},
\]
where \(\operatorname{sat}\) enforces voltage bounds set by the HV driver's output range and the HASEL muscle's dielectric breakdown voltage, and
\(\mathrm{rl}_{\Delta V_{\max}}\{\cdot\}\) is a slew-rate limiter enforcing
\(|u(t)-u(t-T_s)| \le \Delta V_{\max}\) (sampling period \(T_s\)); equivalently
\(|\dot u(t)| \le S_{\max}\) with \(S_{\max}=\Delta V_{\max}/T_s\) (V/s).
We tune \(\Delta V_{\max}\) to limit abrupt command changes that can cause overshoot or disrupt antagonistic clutch scheduling, while respecting hardware slew limits.

\textit{PI feedback:} We employ a PI controller on the tracking error \(e(t)=\theta_d(t)-\theta(t)\):
\[
V_{\mathrm{fb}}(t) = k_p\,e(t) + k_i\int_{0}^{t} e(\tau)\,\mathrm{d}\tau.
\]
Discrete-time realization (sampling \(T_s\)) uses trapezoidal integration for the integral state and conditional integration for anti-windup: the integrator is frozen whenever the tentative command would violate saturation or the slew-rate limit. The measurement is low-pass filtered as \(\tilde\theta_k = \alpha\,\theta_k + (1-\alpha)\,\tilde\theta_{k-1}\) with a small \(\alpha\).

In the feedback-only benchmark, the applied command is
\[
V(t) = \mathrm{sat}_{\underline V,\overline V}\!\big(\mathrm{rl}_{\Delta V_{\max}}\{V_{\mathrm{fb}}(t)\}\big).
\]

\textit{Combined feedforward + feedback:} The combined controller superposes the two contributions before enforcing identical safety constraints:
\[
V_{\mathrm{tot}}(t) = V_{\mathrm{ff}}(t) + V_{\mathrm{fb}}(t),\qquad
V(t) = \mathrm{sat}_{\underline V,\overline V}\!\big(\mathrm{rl}_{\Delta V_{\max}}\{V_{\mathrm{tot}}(t)\}\big).
\]
The feedforward term supplies the bulk of the command from the SSM-based model; the PI term compensates residual model error and disturbances.

\textit{Antagonistic actuation timing with clutches:} The elbow is driven by two antagonistic HASEL muscles, each coupled through an on/off clutch. At any time only one clutch is engaged and receives the command \(V(t)\). We implement a light-weight state machine to schedule side switching safely and smoothly:
\textit{Triggering}: A transition is requested when the reference traverses a dead zone around zero effort (e.g., \(|V_{\mathrm{ff}}|<\varepsilon\) in feedforward/combined modes) or when the desired angle changes sign in feedback-only mode.
\textit{Engagement wait}: Upon trigger, both voltages are set to zero, the target clutch is engaged, and we wait a fixed mechanical engagement time \(\tau_{\mathrm{eng}}\).
\textit{Smooth ramp}: After engagement, the active-side voltage is ramped linearly from 0 to its target magnitude over \(\tau_{\mathrm{ramp}}\) to avoid impulses.
\textit{Cooldown}: A short cooldown \(\tau_{\mathrm{cool}}\) prevents chatter under dithering references.
\textit{Safety enforcement}: All safety checks (saturation, slew-rate, digital outputs in \{0,1\}, and configured joint-angle bounds) are enforced on the final command delivered to hardware.

\textit{Benchmarking protocol:} Each controller variant (feedforward-only, feedback-only, combined) tracks the same reference trajectory under the same safety constraints and clutch scheduler. For fairness and stability, the PI gains were tuned \emph{separately} for the combined and feedback-only conditions (feedforward-only has no gains), targeting a minimal RMS tracking error subject to the same saturation and slew-rate limits.

\backmatter



\bmhead{Acknowledgements}
We are grateful to Joar Axås for the preliminary investigations he carried out in modeling the HASEL actuator from data. We also thank Patricia Apostol for her contribution to the code used to generate the data for the example in section \ref{model_literature}, and Roshan S. Kaundinya for the valuable discussions on applying adiabatic SSM to control.

\bmhead{Author contributions}

L.B., A.K., R.K. and G.H. designed the research. L.B. and A.K. carried out the research. L.B. developed the modeling part and developed the relative software. A.K. performed the experiments and developed the control part. L.B. and A.K. wrote the paper. R.K. and G.H. reviewed the paper. G.H. and R.K. led the research team.

\bmhead{Competing interests}
The authors declare no competing interests.

\bibliography{bibliography}

\end{document}



\title[Nonlinear Spectral Modeling and Control]{%
  Nonlinear Spectral Modeling and Control of\par
  Soft-Robotic Muscles from Data%
}

\subtitle{Supplementary Material}


\author[1]{\fnm{Leonardo} \sur{Bettini}}\equalcont{These authors contributed equally to this work.}

\author[2]{\fnm{Amirhossein} \sur{Kazemipour}}\equalcont{These authors contributed equally to this work.}

\author*[2]{\fnm{Robert K.} \sur{Katzschmann}}

\author*[1]{\fnm{George} \sur{Haller}}

\email{rkk@ethz.ch; georgehaller@ethz.ch}

\affil[1]{Institute for Mechanical Systems, ETH Zurich, Leonhardstrasse 21, 8092 Zurich, Switzerland}

\affil[2]{Soft Robotics Lab, ETH Zurich, Tannenstrasse 3, 8092 Zurich, Switzerland}





\maketitle

\section{General justification for additive slow forcing}\label{additive_method}
Let us consider the non-autonomous dynamical system
\begin{equation}
\left\{
\begin{array}{ll}
\dot{\boldsymbol{x}} =  \boldsymbol{F}\left(\boldsymbol{x}, \boldsymbol{u} \right) & \quad \boldsymbol{F}, \in \mathcal{C}^r, \\[1.2ex]
\dot{\boldsymbol{u}} = \epsilon \boldsymbol{w}\left(\boldsymbol{x}, \epsilon t \right), & \quad 0 \leq \epsilon \ll 1,
\end{array}
\right.
\end{equation}
where $\boldsymbol{x} \in \mathbb{R}^n$ and $\boldsymbol{u}(t) \in \mathbb{R}^f$ denotes a time-dependent external forcing. We use the notation
\begin{equation}
\boldsymbol{F}\left(\boldsymbol{x}, \boldsymbol{0} \right) = \boldsymbol{\tilde{f}}_0\left(\boldsymbol{x}\right).
\end{equation}
Consider now the integral
\begin{equation}
\int_0^1 \frac{d}{ds} \boldsymbol{F}\left(\boldsymbol{x}, \left( 1 - s\right) \boldsymbol{u} \right) ds,
\end{equation}
which can be evaluated as
\begin{equation}\label{LHS_F}
\int_0^1 \frac{d}{ds} \boldsymbol{F}\left(\boldsymbol{x}, \left( 1 - s\right) \boldsymbol{u} \right) ds = \boldsymbol{\tilde{f}}_0\left(\boldsymbol{x}\right) - \boldsymbol{F}\left( \boldsymbol{x}, \boldsymbol{u}\right).
\end{equation}
On the other hand, by the chain rule, we have
\begin{equation}\label{RHS_F}
\int_0^1 \frac{d}{ds} \boldsymbol{F}\left(\boldsymbol{x}, \left( 1 - s\right) \boldsymbol{u} \right) ds = - \left[\int_0^1 \boldsymbol{D}_{\boldsymbol{u}} \boldsymbol{F}\left(\boldsymbol{x}, \left( 1 - s\right) \boldsymbol{u} \right)  ds \right] \boldsymbol{u},
\end{equation}
where $\boldsymbol{D}_{\boldsymbol{u}} \boldsymbol{F}$ denotes the Jacobian of $\boldsymbol{F}$ with respect to the input $\boldsymbol{u}$. Equating \eqref{LHS_F} and \eqref{RHS_F} yields the decomposition
\begin{equation}
\boldsymbol{F}\left(\boldsymbol{x},\boldsymbol{u}\right) = \boldsymbol{\tilde{f}}_0\left(\boldsymbol{x}\right) + \boldsymbol{f}_1\left(\boldsymbol{x},\boldsymbol{u}\right),
\end{equation}
where 
\begin{equation}
\boldsymbol{f}_1\left(\boldsymbol{x},\boldsymbol{u}\right) \doteq \left[\int_0^1 \boldsymbol{D}_{\boldsymbol{u}} \boldsymbol{F}\left(\boldsymbol{x}, \left( 1 - s\right) \boldsymbol{u} \right)  ds \right] \boldsymbol{u}.
\end{equation}
By further splitting the autonomous part into its linear and nonlinear components,
\begin{equation}
\boldsymbol{\tilde{f}}_0\left(\boldsymbol{x}\right) = \boldsymbol{A}\boldsymbol{x} + \boldsymbol{f}_0\left(\boldsymbol{x}\right),
\end{equation}
we recover the form of Eq.~\eqref{nonautonomous_DS} in the main text.

\section{Procedure for the parametrization of the adiabatic Spectral Submanifold and its reduced dynamics from data}\label{data_driven_procedure}

We elaborate here on Section \ref{method_data_driven_aSSM} of the main manuscript, focusing on the data-driven model reduction approach based on the adiabatic Spectral Submanifold (aSSM). Our goal is to obtain a leading-order approximation of the aSSM, $\mathcal{M}\epsilon$. Despite allowing a nonzero $\epsilon$ in Eq. \eqref{nonautonomous_DS} of the main manuscript, we approximate both the aSSM, $\mathcal{M}\epsilon$, and its associated slow manifold (SM), $\mathcal{L}_\epsilon$, with their frozen-limit counterparts $\mathcal{M}_0$ and $\mathcal{L}_0$, which provide close approximations for $\epsilon > 0$ small by the smooth dependence of $\mathcal{M}_\epsilon$ and $\mathcal{L}_\epsilon$ on $\epsilon$ (see \citet{haller_2024}).

To parametrize the aSSM under this approximation, we generate a family of decaying trajectories by perturbing the system around steady states associated with different constant voltage inputs. These perturbations are introduced by applying abrupt voltage jumps relative to the base constant input. Specifically, after allowing the system to settle into a steady state under a fixed voltage, we suddenly change the input to a new constant value. This abrupt change shifts the steady state to a new point in the extended phase space of the $\epsilon=0$ system along the external input $u$ axis, effectively turning the previous steady state into a new initial condition. For each voltage value $u = \bar{u} $, we obtain one or more decaying trajectories, which allow us to extract the local parametrization of the Spectral Submanifold (SSM), $\mathcal{W}\left(E\left(\bar{u}\right)\right)$, and its associated reduced dynamics, as detailed in the Methods \ref{method_data_driven_aSSM} of the manuscript. As illustrated in Fig. \ref{data_driven_procedure}, this procedure is repeated across a range of voltage values, thereby sampling different slices of the extended phase space. The local SSM parametrizations and their corresponding reduced dynamics, obtained at various voltage levels, are then interpolated to construct a global approximation of the aSSM.

\begin{figure*}[!htbp]
    \centering
    \includegraphics[scale = 1]{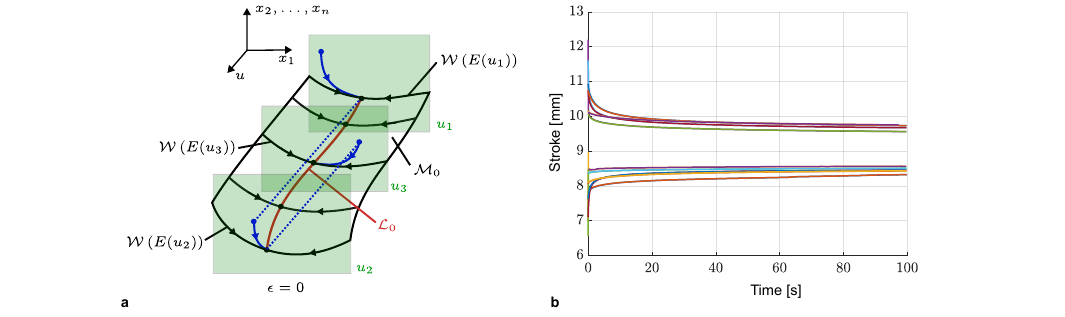}
    \caption{(\textbf{a}) Procedure for the data-driven leading-order parametrization of the aSSM $\mathcal{M}_0$ and the identification of its reduced dynamics. (\textbf{b}) The hysteresis induced by the voltage jumps required to produce distinct initial conditions at a fixed voltage level $\bar{u}$ results in convergence to different fixed points.}
    \label{data_driven_hysteresis}
\end{figure*}

However, the procedure of inducing abrupt voltage jumps to generate initial conditions for decaying trajectories introduces hysteretic effects that are not observed in the normal operational regime of the actuator. This behavior presents a challenge for the practical construction of the aSSM. As seen in Fig. \ref{data_driven_hysteresis}, trajectories initiated from different initial conditions but corresponding to the same reference voltage plane $u = \bar{u} $ in the extended phase space may converge to distinct fixed points. While modeling such behavior becomes relevant in regimes involving extreme voltage changes, it lies beyond the scope of the present work. In our application, these large input jumps are solely employed to populate the data required for aSSM construction and do not reflect standard actuator usage.

Nevertheless, to consistently define the reduced dynamics on the aSSM, it is necessary to associate a unique fixed point with each voltage value $ \bar{u} $. To this end, we calibrate the fixed point selection using the Slow Manifold (SM) obtained from randomly actuated trajectories.

\section{Slow Manifold modeling of antagonistic musculoskeletal joint}

The antagonistic musculoskeletal joint described in Section \ref{data_driven_modeling} of the manuscript employs HASEL actuators in an antagonistic configuration, each coupled with electrostatic clutches, to mimic the actuation of biological limbs (\citet{Kazemipour_2024}). As discussed in the main text, systems of this type exhibit rapid convergence to the SM, which justifies the use of SM-based modeling to simplify the aSSM-based modeling approach. This simplification involves a direct identification of the SM from forced response data as a mapping from input to observable. Specifically, we consider a scalar input $u(t) \in \mathbb{R}$ and a scalar observable $s(t) \in \mathbb{R}$.

A practical and readily interpretable approach is to approximate this mapping using a polynomial in $u$, which is also easily invertible, an advantage for control applications requiring inverse mappings. As an alternative, we have also tested radial basis function (RBF) approximation (\citet{Powell_1992}) with $M = 15$ centers uniformly distributed in the range of the external inputs $[u_{\min}, u_{\max}]$ and variance $\sigma = \left(u_{\max} - u_{\min}\right)/M$. While still interpretable, RBF models can pose challenges in enforcing monotonicity, making their inversion less straightforward.

\begin{figure*}[!htbp]
    \centering
    \includegraphics[scale = 1]{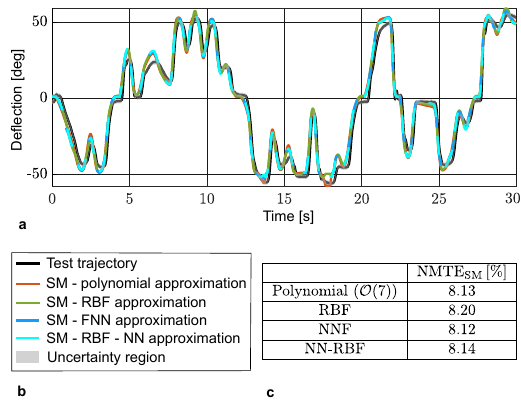}
    \caption{Comparison of different data-driven methods for the one-dimensional SM-based modeling of the antagonistic musculoskeletal system. (\textbf{a}) Time series prediction over a test trajectory. (\textbf{b}) Legend. (\textbf{c}) Normalized Mean Trajectory Errors (NMTEs).}
    \label{data_driven_comparison}
\end{figure*}

We further considered two neural network (NN) architectures: (i) a feedforward neural network (FNN) with one hidden layer of 10 neurons and Tanh activation function (\citet{Hagan_2014}); (ii) a radial basis function neural network (RBF-NN) with a Gaussian activation first layer and linear output layer (\citet{Park_1991}). While these models achieve similar predictions and normalized mean trajectory errors (NMTEs) as the polynomial and the RBF models (see Fig. \ref{data_driven_comparison}), they are less interpretable and more difficult to invert.

\section{Analytical modeling of HASEL actuator}

We consider the HASEL actuator model presented in Section \ref{analytical_models} of the main manuscript. In this appendix, we derive the governing equations for a single pouch of the actuator based on the geometric configuration described in \citet{Kellaris_2019}. Our derivation follows a Lagrangian formalism, which we show to be equivalent to the port-Hamiltonian formalism (\citet{vanderSchaft_2014,Lohmayer_2021}) employed by \citet{Yeh_2022}. 

With reference to Fig. \ref{analytic_geometry_results}(a) of the manuscript, the arc length of the unzipped portion of the pouch is given by
\begin{equation}
l_p(\alpha) = \left( \frac{2A\alpha^2}{\alpha - \sin\alpha \cos\alpha} \right)^{1/2},
\end{equation}
while the zipped length is
\begin{equation}
l_e(\alpha) = L_p - l_p(\alpha).
\end{equation}
The stroke of the actuator, i.e., the contraction from its initial length, is related to $\alpha$ through
\begin{equation}\label{stroke_angle}
x = h - \left( l_p(\alpha) \frac{\sin\alpha}{\alpha} + l_e(\alpha) \right),
\end{equation}
which describes the geometric nonlinearity.
The capacitance of the zipped region is given by
\begin{equation}
C(\alpha) = \frac{\epsilon_r \epsilon_0 w}{2t}, l_e(\alpha),
\end{equation}
where $\epsilon_r$ denotes the relative permittivity of the dielectric film and $\epsilon_0$ is the vacuum permittivity. Using Eq.~\eqref{stroke_angle}, the capacitance can further be expressed as a function of the actuator stroke $x$.

The single pouch is modeled as a single-degree-of-freedom nonlinear oscillator coupled with an electric circuit.
The Lagrangian for the $i$-th pouch is defined as
\begin{equation}
\mathcal{L}_i = T_i - V_i = \frac{1}{2}m \dot{x}_i^2 - \frac{1}{2}k x_i^2 - mg x_i - \frac{1}{2}\frac{Q_i^2}{C\left(x_i\right)}, \qquad i = 1, \dots, N,
\end{equation} 
where $x$ denotes the vertical stroke,  $Q$ is the electric charge accumulated on the electrodes and $C(x_i)$ is the capacitance, which indirectly depends on $\alpha$ through the zipped region. The potential energy includes contributions from the elastic spring, gravity, and the electrostatic energy stored in the capacitor.

To account for dissipative effects, we introduce a Rayleigh dissipation function linear in the vertical velocity, i.e.,
\begin{equation}
\mathcal{R}_i = \frac{1}{2}c \dot{x}_i^2, \qquad i = 1, \dots, N,
\end{equation}
where $c$ denotes the viscous damping coefficient. Hereafter, we omit the subscript $i$ for simplicity.

The equations of motion, given in Eq. \eqref{governing_eq_analytic_HASEL} of the main manuscript, are then obtained as 
\begin{equation}\label{lagrangian}
\frac{d}{dt}\left(\frac{\partial\mathcal{L}}{\partial\dot{x}} \right) - \frac{\partial \mathcal{L}}{\partial x} + \frac{\partial R}{\partial \dot{x}} = 0.
\end{equation}
To complete the coupled electromechanical model, we add the governing equation for the electric charge $Q$ accumulating in the zipped region of the electrodes. This is described by the RC-circuit model 
\begin{equation}
\dot{Q} = -\frac{1}{RC}Q + \frac{1}{R}u.
\end{equation}
driven by an external voltage source $u$. In the literature, the port-Hamiltonian formalism is frequently employed to describe electromechanical systems. Although the term \textit{Hamiltonian} traditionally refers to conservative systems, the port-Hamiltonian approach extends this formalism to include dissipative effects. In this section, we clarify the connection between the dissipative Lagrangian and the port-Hamiltonian formulations by explicitly deriving the equations in a dissipatively perturbed Hamiltonian form.

We define the generalized momentum as
\begin{equation}
p = \frac{\partial \mathcal{L}}{\partial \dot{x}},
\end{equation}
which gives 
\begin{equation}
\dot{p} = \frac{d}{dt}\left(\frac{\partial \mathcal{L}}{\partial \dot{x}} \right).
\end{equation}
From Eq. \eqref{lagrangian}, we obtain
\begin{equation}\label{lagrangian_intermediate}
\dot{p} = \frac{\partial \mathcal{L}}{\partial x} - \frac{\partial \mathcal{R}}{\partial \dot{x}}.
\end{equation}
We now apply the classic Legendre transformation to transition from the Lagrangian to the Hamiltonian setting. This is achieved by expressing $\dot{x}$ as a function of $(x, p)$, i.e., $\dot{x} = \dot{x}(x, p)$, and defining the Hamiltonian function as
\begin{equation}
\mathcal{H}\left(x,p\right) = p \dot{x}\left(x,p\right) - \mathcal{L}\left(x,\dot{x}\left(x,p\right) \right).
\end{equation}
From this definition, we obtain
\begin{align}
\frac{\partial \mathcal{H}}{\partial x} = - \frac{\partial \mathcal{L}}{\partial x},\quad
\frac{\partial \mathcal{H}}{\partial p} = \dot{x}.
\end{align}
Substituting these expressions into Eq. \eqref{lagrangian_intermediate}, we obtain the governing equation for the generalized momentum as
\begin{equation}
\dot{p} = -\frac{\partial \mathcal{H}}{\partial x} - \frac{\partial \mathcal{R}}{\partial \dot{x}}.
\end{equation}
Finally, the coupled evolution equations in a perturbed Hamiltonian form read
\begin{equation}
\begin{array}{l}
\dot{x} = \frac{\partial \mathcal{H}}{\partial p}, \\[2ex]
\dot{p} = -\frac{\partial \mathcal{H}}{\partial x} - \frac{\partial \mathcal{R}}{\partial p}\frac{\partial p}{\partial \dot{x}}.
\end{array}
\end{equation}

\bibliography{bibliography}